\documentclass[11pt]{article}
\usepackage{amssymb, amstex}

\def\mapdownright#1{\Big\downarrow\rlap{$\vcenter{\hbox{$\scriptstyle#1$}}$}}
\def\mapdownleft#1{\rlap{$\vcenter{\hbox{$\scriptstyle#1$}}$}\,\,\Big\downarrow}
\def\mapse#1{\rlap{$\vcenter{\hbox{$\scriptstyle#1$}}$}\,\,\searrow}
\def\mapright#1{\smash{\mathop{\longrightarrow}\limits^{#1}}}
\def\mapleft#1{\smash{\mathop{\longleftarrow}\limits^{#1}}}

\newtheorem{theorem}{Theorem}[section]
\newtheorem{nothing}[theorem]{$\!\!$}
\newtheorem{proposition}[theorem]{Proposition}
\newtheorem{lemma}[theorem]{Lemma}

\newtheorem{construction}[theorem]{Construction}
\newtheorem{remark}[theorem]{Remark}
\newcommand{\prf}{\vspace{.05in}
                    \noindent {\sc Proof} \hspace{.05in}}
\newcommand{\ethrm}{\hspace*{\fill}
                      $\Box$
                      \vspace{.1in}}

\newcommand{\bp}{{\mathbb P}}
\newcommand{\bc}{{\mathbb C}}
\newcommand{\bz}{{\mathbb Z}}
\newcommand{\bq}{{\mathbb Q}}
\newcommand{\br}{{\mathbb R}}
\newcommand{\bd}{{\mathcal D}}
\newcommand{\bl}{{\mathcal L}}
\newcommand{\bn}{{\mathcal N}}
\newcommand{\bh}{{\mathcal H}}

\newcommand{\bigo}{{\mathcal O}}
\newcommand{\bigx}{{\mathcal X}}
\newcommand{\bige}{{\mathcal E}}
\newcommand{\bigw}{{\mathcal W}}
\newcommand{\bigy}{{\mathcal Y}}

\newcommand{\aut}{{\mathop{\rm Aut}\nolimits}}
\newcommand{\cy}{Calabi--Yau}

\newcommand{\cyt}{Calabi--Yau\ threefold}

\newcommand{\spec}{\mathop{\rm Spec\,}\nolimits}
\newcommand{\pic}{\mathop{\rm Pic\,}\nolimits}
\newcommand{\Birat}{\mathop{\rm Birat}\nolimits}
\renewcommand{\phi}{\varphi}
\newcommand{\df}{\mathop{\rm Def\,}\nolimits}
\newcommand{\ext}{\mathop{\rm Ext}\nolimits}
\newcommand{\im}{\mathop{\rm Im}\nolimits}

\begin{document}

\renewcommand{\thefootnote}{$\dagger$}

\begin{center}
{\LARGE \cy\ threefolds with a curve of singularities and} 

\vspace{0.05in}

{\LARGE counterexamples to the Torelli problem}

\vspace{0.2in}

{\large Bal\' azs Szendr\H oi\footnote{Supported by an Eastern European Research Bursary from Trinity College, Cambridge and an ORS Award from the British Government.}}

\vspace{0.15in}

{\large 3 June 1999}

\end{center}

\setcounter{footnote}{0}
\renewcommand{\thefootnote}{\alph{footnote}}
\vspace{-0.1in}
{\small
\begin{center} {\sc abstract} \end{center}
{\leftskip=30pt \rightskip=30pt
Birational Calabi--Yau threefolds in the same deformation family provide a 
`weak' counterexample to the global Torelli problem, as long as they are 
not isomorphic. In this paper, it is shown that deformations
of certain desingularized weighted projective hypersurfaces provide examples
of families containing birational varieties. The constructed examples 
are shown to be non-isomorphic using a specialization argument.\par}

\begin{center}
\end{center}}

\section*{Introduction}

The theory of \cy\ threefolds has received a lot of attention in recent years, 
due to its relation to various fields ranging from birational geometry 
to superstring theory. One of the important problems in the theory is the 

\vspace{0.1in}

\noindent {\bf Torelli Question:} {\it Let $Y_1$, $Y_2$ be smooth, deformation
equivalent \cyt s over $\bc$. Assume that there is an isomorphism
\[(H^3(Y_1), Q_{Y_1}) \cong (H^3(Y_2), Q_{Y_2})\]
respecting the Hodge structures polarized by 
the intersection forms. Are $Y_1$ and $Y_2$ isomorphic?}

\vspace{0.1in}

\noindent The answer to this question remained elusive for a long time 
even in the simplest case, 
for general quintic hypersurfaces in $\bp^4$. However, weak global Torelli in
this case has recently been proved by Voisin~\cite{voisin}. It is possible 
that the answer is positive for simply connected \cy\ threefolds with 
Picard number one. 

In this paper, I give some examples of families of threefolds where the
answer to the above question is negative. They are in a sense `weak' examples.
They arise from a theorem of Koll\'ar, which claims that there
is an isomorphism of polarized  Hodge structures as above 
if $Y_1, Y_2$ are birational \cyt s. In Section~1, I exhibit 
birational threefolds which are at the same time deformation 
equivalent\footnote{For this phenomenon to occur, the Picard number 
must be greater than one, as in the Picard number one case 
the smooth projective birational model is unique.},
based on an idea of~\cite{twoparam} and \cite{enhanced}. 
The varieties arise as deformations of resolutions $Y$ of 
Calabi--Yau threefolds with a curve of singularities. The birational 
isomorphisms give an interesting involution on the base of the Kuranishi 
space of the resolutions $Y$.

In order to obtain a counterexample to Torelli, one has to prove that the 
constructed smooth threefolds are not isomorphic. I will introduce three 
families of explicit examples in Section~2 that exhibit three different sorts
of behaviour. The difference lies in the generic automorphism 
group and the action of this group on the Kuranishi space. In particular,
Proposition~\ref{typec} shows that 
in one family the presence of an involution destroys the counterexample.  
The main result of the paper is Theorem~\ref{maintheorem}, 
which gives explicit cases where (weak) global Torelli fails. 
In particular, a question posed in~\cite[1.6]{rolling} is answered. 
The proof is based on a specialization argument and a standard
result on automorphisms in families.

The construction shows that one can find threefolds with the same 
Hodge structure in an arbitrary small disc neighbourhood of the central
fibre in the Kuranishi space of $Y$. This however does not contradict 
Infinitesimal Torelli: as explained in Remark~\ref{infto}, the period point
of $Y$ in the period domain is fixed by an element of the 
arithmetic monodromy group.

Desingularizations of Horrocks--Mumford quintics give another example
of birational, deformation equivalent \cy\ threefolds~\cite{borcea_hm},  
the Picard number being four in that case. Indeed, it is reasonable to 
expect that as the Picard number increases, Torelli can fail more and more 
badly, although under a weak condition there will always be
finitely many isomorphism classes with the same Hodge structure 
within a deformation family by the main result of~\cite{finite}. 

The examples suggest the 

\vspace{0.1in}
\noindent {\bf Modified Torelli Question:} 
{\it Let $Y_1$, $Y_2$ be smooth, deformation
equivalent \cyt s over $\bc$, with isomorphic polarized Hodge structures. 
Are $Y_1$ and $Y_2$ birationally equivalent?}
\vspace{0.1in}

The example of Aspinwall--Morrison~\cite{am}, featuring non-simply connected 
threefolds, may provide a counterexample to this more general question, 
but some details of that example remain to be worked out\footnote{A detailed
analysis of the example will be given in~\cite{thesis}.}. 

My final remark is that the theory of holomorphic symplectic varieties
shows somewhat parallel features. On one hand, there is an example of 
birational, non-isomorphic varieties~\cite{debarre}. On the other hand,
the polarized Hodge structures (on second
cohomology) of birational varieties are isomorphic, and the natural
way to pose the Torelli problem is exactly the same as above~\cite{huy}. 
However, similarities end here, \cite[4.6]{huy} shows that two birational 
holomorphic symplectic manifolds can always be realized as central fibres
in families over the disc, the families being isomorphic over the 
punctured disc. This phenomenon is special to the holomorphic symplectic case.

\vspace{0.1in}

\noindent
{\bf Acknowledgements} \ I wish to thank Pelham Wilson for suggesting the idea
of the construction, and also for encouragement, comments and corrections.
Thanks are also due to Miles Reid 
for his interest in my work and friendly advice,
to Michael Fryers, Daniel Huybrechts, Umar Salam and Burt Totaro 
for discussions and comments and finally to the referee for 
suggesting some improvements and pointing out Proposition~\ref{topology}. 

\section*{Notation and definitions} 

A {\it \cy\ threefold} is a normal projective threefold $Y$ with canonical 
Gorenstein singularities, satisfying $K_Y\sim 0$ and $H^1(Y,\bigo_Y)=0$. 
The {\it nef cone} of $Y$ is the closed cone generated by ample classes in
$\pic(Y)\otimes \br \cong H^2(Y,\br)$. 

The {\it weighted projective space} $\bp^m[w_0,\ldots w_m]$ is 
the quotient of $\bc^{n+1}\setminus\{0\}$ 
by the $\bc^*$-action having the given 
weights. A {\it weighted hypersurface}  
$X=X_d\subset\bp^m[w_0,\ldots w_m]$ is defined by the 
vanishing of a homogeneous polynomial $f$ of weighted degree $d$.

If $Y$ is a variety, by a slight abuse
of notation I write $T_0\df(Y)$ for the linear space classifying 
first order deformations of $Y$. This space is isomorphic to 
$H^1(Y,\Theta_Y)$. 
If $f:Y\rightarrow X$ is a morphism of algebraic varieties, a
{\it deformation of $f$ over $S$} is a commutative diagram 
\[
\begin{array}{ccccccc}
Y &&&\longrightarrow&&&\bigy\\
&\mapse{f}&&&&\swarrow \\
\Big\downarrow&&X&\rightarrow & \bigx & &\mapdownright{g}\\
& \swarrow &&&&\mapse{h} \\
0 &&& \hookrightarrow &&& S
\end{array}
\]
where $g,h$ are flat and $Y=g^{-1}(0)$, $X=h^{-1}(0)$. 
There is a vector space $T_{0}\df(Y,f,X)$ classifying the first
order deformations, defined defined by Ran~\cite{ran} 
as a complicated $\ext$ group over a Grothendieck
topology. There are natural maps
\[ T_{0}\df(X)\leftarrow T_{0}\df(Y,f,X)\rightarrow T_{0}\df(Y)\]
to the first order deformation spaces of $X$ and $Y$. 

Finally, $H^n(X)$ denotes the integral cohomology of $X$ modulo torsion.

\section{The general framework}

The counterexamples to the Torelli problem presented in this paper
will be based on the following result:

\begin{theorem} {\rm (Koll\'ar~\cite[4.12-13]{flops}.)} 
Assume that $Y_1, Y_2$ are smooth \cy\  threefolds,
related by a flop $\psi: Y_1 \dashrightarrow Y_2$. 
Then $\psi$ induces an isomorphism of polarized Hodge structures
\[(H^3(Y_1), Q_{Y_1}) \cong (H^3(Y_2),Q_{Y_2}).\]
\label{kollar_thrm}
\end{theorem} 
As a first step, I introduce a class of \cyt s having birational 
varieties in their deformation space.

\begin{construction} \rm Let $Y$ be a smooth \cyt\ containing 
a surface $E$ ruled over a curve $C$ of genus $g\geq 2$, necessarily smooth 
by \cite{wilson},~\cite{wilson_symplectic}, which 
can be contracted inside $Y$ by a log-extremal contraction given by the 
divisor $H\in\pic(Y)$:

\[
\begin{array}{rccl}
\phi_{\mid mH \mid}:& Y & \rightarrow & X \\
&\cup && \cup \\
& E & \rightarrow & C.
\end{array}
\]
For the rest of the paper, let
$\pi:\bigy\rightarrow S$ be the the Kuranishi family of $Y$. By
Unobstructedness \cite{tian},~\cite{todorov} its base $S$ is smooth.  
Fix once and for all an identification between $S$ and an
open disc in $T_0\df(Y)$; this defines a notion of 
(partial) addition and scalar multiplication on $S$. 
There is an obvious one-to-one correspondence 
between vectors $v\in S$ 
and maps from the unit disc $f_v:\Delta\rightarrow S$, given by 
$f_v(z)=zv$. The map $f_v$ gives a 
pullback family $\pi_v:\bigy_v\rightarrow \Delta$ with smooth \cy\ fibres.

Let
\[S_E = \im\left(T_0\df(E,k,Y)\rightarrow T_0\df(Y)\right)\cap S,\]
this is a linear subspace of $S$, corresponding to the deformation 
directions along which $E$ deforms together with the deformation.
Wilson \cite{wilson} shows that if $E$ is smooth, 
$S_E\subset S$ is of codimension~$g$. 

Let $\pi_v:\bigy_v\rightarrow \Delta$ be a general one-parameter deformation
for $v$ not contained in $S_E$. Assume that the nef cone is invariant 
in the family. (This will be proved for all the discussed examples 
in Lemma~\ref{cone}.) The following is well-known:
\end{construction}

\begin{proposition} There exists a unique relative Cartier divisor $\bh$ 
on $\bigy_v$ extending $H$ on $Y$. There is a morphism  
$\phi_{\mid m\bh \mid}:\bigy_v\rightarrow\bigx_v$ over $\Delta$,
which is the contraction of
$E$ in the central fibre and the contraction of a finite number 
of rational curves $\cup_i C^i_t$ on the general fibre $Y_t$, 
the number of such curves being $2g-2$ if counted with 
appropriate multiplicities. The map $\bigx_v\rightarrow\Delta$ is flat.
\end{proposition}
\prf As $h^2(\bigo_Y)=0$, the divisor $H$ extends uniquely over the family, 
and by assumption $\bh_t$ is a nef and big divisor on all fibres $Y_t$. 
So $\pi_{v*}(m\bh)$ is a vector bundle over $\Delta$ and 
$\mid m\bh \mid$ defines a morphism. The exceptional locus over 
$t\in\Delta^*$ is a finite union of rational curves by \cite{wilson}, 
their expected number (Gromov-Witten invariant) is calculated in
\cite[2.3]{wilson_flop}.
The last statement follows from \cite[11.4]{kollar-mori}. 
\ethrm

\begin{proposition} {\rm (cf.~\cite[2.3]{wilson_flop}.)} 
There exists a flop of $\bigy_v\rightarrow \Delta$, i.e. a diagram
\[
\begin{array}{ccccc}
\bigy_v & & \stackrel{\psi}{\dashrightarrow } && \bigy^+_v \\
&\searrow & & \swarrow & \\
&& \bigx_v  
\end{array}
\]
over $\Delta$, where the birational map $\psi_t$ flops 
the curves $C_t^i$ on $Y_t$ for $t\in\Delta^*$ and gives 
an isomorphism on the fibre over $0$.
\end{proposition}
\prf $\bigx_v$ has only cDV singularities, so the existence of the flop can be
seen by taking hyperplane sections and using \cite[11.10]{kollar-mori}.  
$\bigy^+_v$ is smooth, so the birational map 
$\bigy_v \dashrightarrow \bigy^+_v$ 
is an isomorphism in codimension 1, which shows that $\psi_0$ is an 
isomorphism.
\ethrm

The family $\bigy^+_v\rightarrow \Delta$ has central fibre $Y$, so corresponds
to some vector $\alpha (v) \in H^1(Y,T_Y)$. 
Shrinking $S$ if necessary and defining $\alpha$ as the identity on $S_E$, 
one obtains a map $\alpha: S\rightarrow S$ which is clearly an 
involution. The following is a tautology: 

\begin{proposition} The map $\alpha$ is linear in the obvious sense. Its 
fixed locus is exactly $S_E$. 
\label{fixlocus}
\end{proposition} 
\prf Let $v_1, v_2\in S\setminus S_E$ such that $v=v_1+v_2\in S$ also, 
then the maps
$f_{v_i}$ together with the linear structure on $S$ 
give a map $\Delta^2\rightarrow S$. The map $f_v$ coincides
with the composite map $\Delta\rightarrow\Delta^2\rightarrow S$ where the
first map is the diagonal one. The flop can just as well be constructed over
$\Delta^2$ which shows linearity immediately. The obvious modification of
this argument works also if some $v_i\in S_E$. 

Now suppose $v$ is not in $S_E$ but it is fixed by $\alpha$. By universality 
of the Kuranishi family, this means that the birational map $\psi$ is 
the identity on all fibres. This is clearly nonsense.
\ethrm

\section{Particular families}

In this Section, I will briefly investigate the geometry of 
three families which fit 
into the framework described above. The singular variety $X$ will be 
a general hypersurface $X_8\subset\bp^4[1^2,2^3]$, 
$X_{12}\subset\bp^4[1^2,2^2,6]$ or $X_{14}\subset\bp^4[1,2^3,7]$, 
respectively. $X$ has Picard number one, canonical 
singularities along curves and by adjunction, trivial canonical sheaf. 

\begin{nothing}\label{firstcase} \rm
First assume that $X = X_8 \subset \bp^4[1^2,2^3]$,
a variety discussed at length from the point of view of
mirror symmetry in \cite{twoparam}, deformations of which also featured
in~\cite[1.6]{rolling}. $X$ is singular along the locus $C=\{x_i=0\}$, a plane
curve of genus 3. 

The linear system $\bigo(2)$ embeds $\bp^4[1^2,2^3]$ as a quadric 
$Q_3=\{z_1z_3=z_2^2\}$ of rank $3$ in $\bp^5$, 
where $z_i$ are the coordinates on the $\bp^5$. 
The hypersurface $X$ becomes 
a complete intersection of $Q_3$ with a quartic $F_4$ in $\bp^5$.
The singularities of the quadric $Q_3$ can be resolved
by the map $F(2,0,0,0)\rightarrow Q_3$,  
where $F(2,0,0,0)$ is a rational scroll. $F(2,0,0,0)$ generically deforms to
the scroll $F(1,1,0,0)$. The contraction 
$\varphi:F(1,1,0,0)\rightarrow Q_4\subset\bp^5$ 
maps the scroll to a quadric $Q_4=\{z_1z_3=z_2^2-tz_4^2\}$ of rank 4. 
Intersections of $Q_4$ with $F_4$ give deformations $X_t$ of $X$, 
with four isolated $cA_1$ points at $\{z_1=\ldots =z_4=0\}$.
The resolution $Y_t\rightarrow X_t$ replaces these points by 
the rational curves $C^i_t$. 
The following diagram summarizes the state of affairs:
\[
\begin{array}{ccccccccccccc}
&&F(1,1,0,0) &\supset &Y_t & \leadsto & Y & \subset & F(2,0,0,0) \\
&&\downarrow & &\downarrow&&\downarrow && \downarrow \\
\bp^5 & \supset&Q_4 & \supset & X_t & \leadsto & X & \subset & \bp^4[1^2,2^3]& \cong & Q_3 & \subset \bp^5
\end{array}
\]

$F(1,1,0,0)$ is a quotient of 
$(\bc^2\setminus\{0\})\times(\bc^4\setminus\{0\})$ by the group
$(\bc^*)^2$, where the two multiplicative actions have weights
$(1,1;-1,-1,0,0)$ and $(0,0;1,1,1,1)$ respectively. If $t_i$,
$u_i$ denote the coordinates on the affine spaces, 
then the map $\varphi$ is given by
\[
\varphi: (t_1,t_2;u_1,\ldots, u_4)\mapsto (u_1t_1:\frac{1}{2}(u_1t_2+u_2t_1):u_2t_2:\frac{1}{2\sqrt{t}}(u_1t_2-u_2t_1):u_3:u_4),
\] 
having the quadric $Q_4$ as its image. 
There is however an ambiguity in the choice of the sign of the square 
root. This does not matter in $\bp^5$, as the two choices are isomorphic under
the map $\sigma:\bp^5\rightarrow\bp^5$, 
$z_4\mapsto -z_4$, $z_1\leftrightarrow z_3$. 
However, taking the intersection with the
quartic $F_4$, the two choices give different resolutions as long as 
$\sigma$ is not an automorphism of $F_4$. This gives the
two resolutions related by a flop.

Note that if $\sigma$ is an automorphism of $F_4$, then it gives an 
automorphism $j\in\aut(Y)$ and the map $\alpha$ acting on the base of the
Kuranishi space $S$ of $Y$ equals $j^*$. For such special points in the moduli
space, the families $\bigy_v$, $\bigy_{\alpha(v)}$ are isomorphic. 
\end{nothing}

\begin{nothing}\label{secondcase} \rm
Next consider the other two families $X_{12}\subset\bp^4[1^2,2^2,6]$ and
$X_{14}\subset\bp^4[1,2^3,7]$. These varieties are singular along the curves 
given by the vanishing of the variables of odd degrees, of genus 
2, 15 respectively. 

Notice that these varieties always
have nontrivial automorphisms. If $z$ denotes the variable of highest 
degree and $x_i, y_i$ the other variables, then changing variables
$X=\{z^2+f_d(x_i,y_i)=0\}$ and then $z\mapsto -z$ gives an
involution $i$. 
The automorphism $i$ of $X$ extends to an involution $j$ on the 
resolution $Y$. $j$ induces a natural linear action $j^*$ 
on the space $T_0\df(Y)$; denote its fixed locus
by $T^+_0\df(Y)$. This is the subspace of $T_0\df(Y)$ 
corresponding to deformation directions 
along which the involution $j$ also deforms.
Notice that by universality of the Kuranishi family, 
the families $\bigy_{v}\rightarrow \Delta$, 
$\bigy_{j^*(v)}\rightarrow \Delta$ are isomorphic under an isomorphism 
induced by $j$. The following is again a tautology:

\begin{proposition} The actions of $\alpha$, $j^*$ on $S$ commute.
\end{proposition}

\prf The following diagram obviously commutes:
\[
\begin{array}{ccccc}
\bigy_v &  \stackrel{\psi}{\dashrightarrow } & \bigy^+_v \\
\Big\downarrow& &\Big\downarrow  \\
\bigy_{j^*(v)} &  \stackrel{\psi}{\dashrightarrow } & \bigy^+_{j^*(v)}. 
\end{array}
\]
Thus $\bigy_{j^*\alpha(v)}=\bigy^+_{j^*(v)}=\bigy_{\alpha j^*(v)}$.
\ethrm

In the two cases, the involution $j^*$ behaves very differently: 

\begin{proposition} 
If $Y$ is the resolution of $X_{12}\subset \bp^4[1^2,2^2,6]$, 
then $j^*$ is trivial. If $Y$ is the resolution of 
$X_{14}\subset \bp^4[1,2^3,7]$, then $j^*=\alpha$ as maps acting on the
base of the Kuranishi space of $Y$. 
\label{typec}  
\end{proposition}

\prf Let $X =X_{12}\subset \bp^4[1^2,2^2,6]$, a double cover 
$X\rightarrow \bp^3[1^2,2^2] \cong Q_3\subset\bp^4$
of a quadric of rank 3 in $\bp^4$, branched over a sextic. 
The quadric is resolved by the scroll $F(2,0,0)\rightarrow Q_3$, 
and $Y$ is a branched double cover of this scroll. The 
variety $Y$ appears in a paper of Fujita~\cite{fujita}, 
where it is denoted by $\Sigma^3(2,0,0)^+_{3,0}$. 
By~\cite[7.13]{fujita}, the action of $j^*$ on $S$ is trivial.
This proves the first statement. 

Now let $X=X_{14}\subset\bp^4[1,2^3,7]$. The fixed locus 
of the involution $i:z\rightarrow (-z)$ is in this case
reducible; the quotient is
$X/\langle i\rangle \cong \bp^3[1,2^3]\cong\bp^3$. The image of the curve
$C$ is $\{t_1=0\}\cap \{g_7(t_1,\ldots t_4)=0\} = \gamma$.
The involution extends to the resolution $Y$ as an involution~$j$. 
The quotient upstairs is $Y/\langle j\rangle\cong W=B_\gamma\bp^3$, 
the blowup of $\bp^3$ along the curve $\gamma$. 
One obtains a diagram 
\[
\begin{array}{ccccc}
E & \longrightarrow &  E^\prime &&\\
\mapdownleft{k} & & \mapdownright{l} &&\\
Y &  \longrightarrow &  W &=& B_\gamma \bp^3 \\
\mapdownleft{f} & & \mapdownright{\pi}  &&\\
X &  \longrightarrow &  \bp^3 && \\
\cup && \cup \\
C & \longrightarrow & \gamma.  
\end{array}
\]

Let $\Theta_W$ be the tangent bundle of $W$, $N_{\gamma/\bp^3}$ the normal
bundle of $\gamma$ in its ambient space. Standard arguments show

\begin{lemma} The natural map 
$H^0(\gamma,N_{\gamma/{\bp^3}})\rightarrow H^1(W,\Theta_W)$ 
is surjective, i.e. any (first-order)  
deformation of $W$ comes from a deformation of $\gamma$ in $\bp^3$. 
\label{blowup}
\end{lemma} 
\ethrm

\begin{lemma} There is an inclusion
\[T^+_0\df(Y) \subset  \im\left(T_0\df(E,k,Y)\rightarrow T_0\df(Y)\right).\] 
\label{inclusion}
\end{lemma}
\prf Consider
\[
\begin{array}{ccccc}
T_0\df(E,k,Y)& \mapleft{\lambda} &  F & \longrightarrow & T_0\df(E^\prime, l,W)\\
\mapdownleft{\nu} && \mapdownleft{\beta} && \mapdownright{\delta}\\
T_0\df(Y) & \mapleft{\psi} &T_0^+\df(Y) & \mapright{\mu}& T_0\df(W).
\end{array}
\]

\noindent Here $\mu$, $\delta$, $\nu$ 
are the obvious maps and $\psi$ is the inclusion. Let   
$F$ be the fibre product; its elements are pairs of vectors giving a
first-order 
deformation $\bigy\rightarrow B=\spec \bc\left[\epsilon\right]/(\epsilon^2)$ 
of $Y$ with an involution $J$ on $\bigy$, 
and a first-order deformation 
$(\bigw, \bige^\prime)\rightarrow B$ of the inclusion
$l:E^\prime\rightarrow W$, with compatible image in $T_0\df(W)$. 
In other words, elements of $F$ give a diagram 
\[
\begin{array}{ccccc} 
&& \bige^\prime\\
&&\Big\downarrow\\
\bigy & \mapright{q} & \bigy/\langle J \rangle &= & \bigw \\ 
& \searrow & \Big\downarrow  \\
&& B
\end{array}
\]

\noindent where $q$ is the quotient map and all maps over $B$ are flat. 
$q$ is finite, so $\bige=q^{-1}(\bige^\prime)$
is a relative Cartier divisor on $Y$, flat over $B$, hence one
obtains a deformation $\bige\rightarrow\bigy$ of 
$E\rightarrow Y$. This defines the map $\lambda$ in the previous diagram.
By construction, the left square becomes commutative. 

Lemma~\ref{blowup} shows that $\delta$ is surjective, so $\beta$ must also 
be surjective.
Then the commutativity of the left square shows that the image of the 
inclusion $\psi$ must be contained in the image of $\nu$. This proves 
Lemma~\ref{inclusion}. 
\ethrm

\noindent To complete the proof of Proposition~\ref{typec}, notice that 
$j^*=\alpha$ would follow from the equality of fixed loci, 
as the two maps are commuting involutions, linear with respect to the
partial linear structure on~$S$. 
The fixed locus of $\alpha$ is $S_E$, whereas the fixed locus of
$j^*$ is by definition $T_0^+\df(Y)\cap S$. 
By the previous Lemma, $S_E\supset T_0^+\df(Y)\cap S$ and
it suffices to show that the dimensions here are equal. The dimension of 
$S_E$ is $h^{2,1}(Y)-g(C)=107$. On the other hand, by a standard dimension 
count, the 
hypersurface $X_{14}\subset \bp^4[1,2^3,7]$ depends on $107$ parameters.
Resolutions of such hypersurfaces are always double covers, so 
$\dim T_0^+\df(Y)\geq 107$. This concludes the proof. 
\ethrm
\end{nothing}

Thus in the last case, the varieties $Y_t$, $Y_{\alpha(t)}$ are isomorphic
for all $t$. For the rest of the paper, restrict attention to the first two 
families. The last statement in this Section describes 
the nef cone of $Y$ in these cases; this description will be needed below. 

\begin{proposition} The nef cone of the resolution $Y$ is
generated by the flopping face and the face corresponding to a fibration with
base $\bp^1$. It is constant on any deformation of $Y$.   
\label{cone}
\end{proposition}
\prf The first statement is clear from the above discussion; the map to $\bp^1$
comes from the structural map of the rational scrolls. 
As for the second part, by the main result of~\cite{wilson}, 
the nef cone is invariant on all deformations if no deformation of $Y$ 
contains a surface which is (quasi-) ruled over an elliptic curve. 
A sufficient condition 
for this is that there is no class $F\in H^2(Y,\bq)$ satisfying 
$F^3=c_2(Y)\cdot F=0$. $H, E\in H^2(Y,\bq)$ are not linearly dependent, so
they give a $\bq$-basis of the rank-two space $H^2(Y,\bq)$ and the
existence of the class $F$ is equivalent to 
$((c_2(Y)\cdot H)E-(c_2(Y)\cdot E)H)^3=0$.  A routine calculation
shows that this fails in both cases.
\ethrm

\begin{remark}\rm The families investigated in this Section belong to 
a larger set of examples that can be found using a systematic search
based on the theory of graded rings of weighted complete intersection 
varieties. For details, consult~\cite{thesis}. 
\end{remark}

\section{The automorphism group of the general variety}

\begin{proposition} The general
$X_8 \subset \bp^4[1^2,2^3]$ has trivial automorphism group, 
whereas the general
$X_{12} \subset \bp^4[1^2,2^2,6]$ has automorphism group $\bz/2\bz$. 
\label{autom}
\end{proposition}

\prf Assume that $\sigma$ is a 
nontrivial element of the automorphism group of the hypersurface 
$X=\{F_8(\vec x, \vec y)=0\}\subset \bp^4[1^2,2^3]$,
\[
F_8 (\vec x, \vec y) =  f_8(\vec x) + g_1 (\vec y) f_6(\vec x) + g_2 (\vec y) f_4(\vec x) + g_3 (\vec y) (x_1^2+x_2^2)+g_4 (\vec y),
\]
where $\vec x=(x_j), \vec y=(y_j)$ are the homogeneous coordinates 
of degrees $1,2$ respectively.  
The Picard group of $X$ is of rank one, which implies
that $\sigma$ comes from a projective automorphism of 
$\bp^4[1^2,2^3]$. Hence by \cite[4.7]{cox}, 
\begin{eqnarray*}
\sigma(\vec x)& = & A\vec x, \\
\sigma(\vec y)& = & B\vec y + C (S^2\vec x),
\end{eqnarray*}
where $A$ is a $2\times 2$ matrix, $B, C$ are $3\times 3$ matrices and 
$S^2\vec x$ is a shorthand for $(x_1^2,x_1x_2,x_2^2)^t$. 

The singular locus of $X$ is 
the genus three plane curve $\{g_4=0\}\subset \bp^2$. 
This has to be mapped isomorphically by $\sigma$, so if $g_4$ is general, then
after fixing an overall constant, $B=I$.
Writing out the conditions for the invariance of the
cubic $y_i$ terms, one obtains $C=0$ for general choice of $F$.  
Finally, the automorphism 
has to fix the octic $f_8(\vec x)=0$, so for general $f_8$, 
the only possibility is $A=\delta I$ and then clearly
$\delta=\pm 1$. Thus up to constant, $A=\pm I, B=I, C=0$, where the final sign 
is part of the $\bc^*$-action in the definition of weighted projective space.

The proof for a general $X_{12} \subset \bp^4[1^2,2^2,6]$ is completely 
analogous, so it is omitted. 
\ethrm

\section{The conclusion}

I will need the following rather standard result: 

\begin{proposition} Let $\bigx\rightarrow B$ be a family of \cy\ threefolds
with canonical singularities over a complex space $B$, 
having a simultaneous resolution
$\bigy\rightarrow\bigx$ over $B$. Let $\bl$ be a relatively ample
relative Cartier divisor on $\bigx$. Let $\aut_B(\bigx, \bl)$ 
be the scheme of relative automorphisms in the family. 
Then $\aut_B(\bigx, \bl)$ is finite and unramified 
over $B$. \label{fp}
\end{proposition}
\prf $\aut_B(\bigx, \bl)$ is unramified over $B$, 
as the fibres of the family $\bigx\rightarrow B$ are 
varieties without infinitesimal automorphisms. Quasi-finiteness is clear, and
properness follows from the valuative criterion 
(see e.g.~\cite[Lemma 4.2]{fantechi-pardini}).
\ethrm

\begin{proposition} Let $X$ be a general member of one of the 
two families, $Y$ the \cy\ resolution. Let $v\in S\setminus S_E$ be 
a deformation direction as in Section~1, and $\bigy\rightarrow \Delta$
the corresponding family with flop $\bigy^+\rightarrow \Delta$.
Assume that for all $t\in U$ in a dense set $U\subset \Delta^*$, 
there exists an isomorphism $Y_t\cong Y^+_t$. Then the variety $X$ 
has nontrivial automorphism group, respectively 
automorphism group larger than ${\bz/2\bz}$.
\label{reduce}
\end{proposition} 

\begin{lemma} Assume that for some $t\in\Delta^*$, there exists an
isomorphism $Y_t\cong Y^+_t$. Then $\aut(X_t)$ is nontrivial, 
respectively larger than $\bz/{2\bz}$.
\end{lemma}
 
\prf Assume that $Y_t\cong Y^+_t$, then the flop 
corresponds to a nontrivial birational automorphism $\psi\in\Birat(Y_t)$, 
not the identity on the complement of the exceptional
locus (and not the involution $j$).
So it descends to a nontrivial birational automorphism 
$\bar\psi\in\Birat(X_t)$. On the other hand, using the fact that
any isomorphism $Y_t\cong Y_t^+$ must identify faces of the nef cones 
of the same type, it is easy to check that 
$\bar\psi$ must fix a suitable multiple of the ample 
generator of the Picard group of $X_t$. Hence $\bar\psi$ is biregular by
\cite[2.1.6]{kollar_etc}.
\ethrm

\noindent {\sc Proof of Proposition \ref{reduce}} \hspace{.05in}
Let $\bigx\rightarrow \Delta$ be the contracted family, 
$\bn$ a relatively ample sheaf on $\bigx$.  
If there exists an isomorphism $Y_t\cong Y^+_t$ for all
$t\in U$, there is a nontrivial element in  
$\aut_U(\bigx,\;\bn)$ for all $t\in U$. The statement now follows 
from Proposition~\ref{fp}.  
\ethrm

Let $\bd$ denote the period domain parameterizing
polarized Hodge structures on the $\bz$-module $H^3(Y)$. Let 
\[\Gamma=\aut(H^3(Y), Q_Y)\] 
be the corresponding arithmetic monodromy group. 
The following is the main result of the paper:

\begin{theorem} Let $Y$ be a resolution of a general 
$X_8 \subset \bp^4[1^2,2^3]$ or $X_{12} \subset \bp^4[1^2,2^2,6]$. 
Then (weak) global Torelli fails for $Y$:
the period map 
is finite of degree at least two from the deformation space of 
$Y$ modulo isomorphisms onto its image in $\bd/\Gamma$, the period domain 
modulo monodromy.   
\label{maintheorem}
\end{theorem}
\prf Finiteness of the period map follows from~\cite[4.3]{finite}, 
keeping in mind Lemma~\ref{cone} above. 
The rest follows from Theorem~\ref{kollar_thrm}, 
Proposition~\ref{autom} and Proposition~\ref{reduce}.
\ethrm

\begin{remark} \rm Fixing a marking of the cohomology of $Y$ and
using the Gauss-Manin connection in the bundle $R^3\pi_*\bc$ over 
the base of the Kuranishi space $S$, one obtains the (local) 
period map $\phi:S\rightarrow \bd$. By the Infinitesimal Torelli 
theorem, $\phi$ is an embedding. However, this does not contradict the above
result: the period map is not invariant under the involution
$\alpha$ acting on $S$, 
it is only equivariant with respect to an element $\gamma\in\Gamma$ 
fixing the period point of the central fibre $Y$.
To conclude the paper, I give a geometric description of $\gamma$. 

There is a split exact sequence of Hodge structures 
\[ 0\rightarrow H^3(X) \rightarrow H^3(Y) \rightarrow H^1(C)[-1] \rightarrow 0\]
where the first map is pullback, whereas the second is the dual of the 
cylinder homomorphism given by the family $E\rightarrow C$ of rational
curves in $Y$. Note that the Hodge structure on 
$H^3(X)$ is pure, as $X$ has only quotient singularities. 
\label{infto}
\end{remark}

\begin{proposition} The element $\gamma\in\Gamma$ is 
the involution of the Hodge structure $H^3(Y)$ that fixes $H^3(X)$ and
reflects the sub-Hodge structure $H^1(C)[-1]$ generated by the
family $E\rightarrow C$ of rational curves.
\label{topology}
\end{proposition}
\prf By construction, the image of
$H^3(X)$ is fixed by the involution $\gamma$. On the other hand, 
it is clear that $\gamma$ induces the same action as $\alpha$  
on the space $H^{2,1}(Y)\cong H^1(Y,T_Y)$, the latter isomorphism 
being well-defined up to a constant. Thus the corank of
the submodule fixed by $\gamma$ is at least $2g$ by Lemma~\ref{fixlocus}. 
Hence the fixed submodule is exactly $H^3(X)$ and this concludes the proof.
\ethrm

\noindent {\small \sc Department of Pure Mathematics and Mathematical Statistics

\noindent University of Cambridge

\noindent 16 Mill Lane, Cambridge, CB2 1SB, UK 

\noindent Fax: +44 1223 337 920

\noindent E-mail address: \tt balazs@@dpmms.cam.ac.uk}

\end{document}